# Systems of energy emitting bodies and their properties

## Vitaly O. Groppen

(*North Caucasian Institute of Mining and Metallurgy, Vladikavkaz, Russia*)

Proposed is system of consistent mathematical models describing physical laws of a system of energy emitting bodies in dynamics, relativity and nuclear physics. It is shown the use of developed models for the description of systems, consisting of stable as well as of radioactive bodies and permitting to improve the quality of predicting the binding energy of light stable nuclei using modified semi-empirical equation. Experimental verification of proposed approach with respect to some nuclei of the Periodic Table elements in the first approximation confirms its validity, while we cannot as yet consider it statistically significant.

1. **Introduction**

The present work aims at the analysis of the properties of a system of energy emitting bodies independent of this energy nature. The approach described is based on systems being stable from the point of view of an observer located within such systems, consisting of components which are unstable from the viewpoint of an exterior observer, i.e. consisting of energy emitting bodies. At the same time:
- system of bodies is assumed to be stable if from the point of view of an interior observer the system satisfies the conservation laws;
- body is assumed to be unstable if an exterior observer confirms that the mass (energy) of this body is a function of time.

In other words the author aims at
- creation and investigation of a system of consistent mathematical models which should describe physical laws applicable to the systems consisting of material particles emitting energy;
- use of developed models for the description of real physical processes.

2. **Basic conceptions**

The following assumptions are believed to be valid with respect to the studied below models:
1. Energy $E_B$ concentration in a certain material point B is inseparably linked with its emission to the closest area being in the direct proximity to B, where energy concentration is lower.
2. The emission intensity in the direction of $\varphi$ is directly proportional to the difference of energy concentrations in point, corresponding to particle B and in its immediate vicinity in the direction of $\varphi$.
3. Any polynomial functions of n variables used for description of physical processes (n>0), $f(x_1, x_2, ..., x_n)$, for short time intervals can be described by a polynomial, with integer and non-negative degrees of its variables not exceeding "1":

$$\lim_{t \to 0} f(x_1, x_2, \ldots, x_n) \approx \sum_{i_1=0}^{1} \sum_{i_2=0}^{1} \ldots \sum_{i_n=0}^{1} c_{i_1 i_2 \ldots i_n} \cdot \prod_{j=1}^{n} x_j^{i_j}; \qquad (1)$$

Denoting the surface of the body corresponding to B material particle by symbol "s" the energy loss intensity by this body can be determined by the following:

$$\oint \frac{\partial^2 E_B}{\partial s \partial t} ds = \frac{dE_B}{dt}, \qquad (2)$$

the second member of (2) in accordance with the above assumptions in the isotropic space can be described by the polynomial below:

$$\frac{dE_B}{dt} = k_0 - k_1(E_B - E_\varphi), \qquad (3)$$

where $k_0$, $k_1$ are the coefficients; $E_\varphi$ is energy density in the immediate vicinity of B in the direction of $\varphi$.

Assuming $E_B \gg E_\varphi$, (3) can be transferred to $\frac{dE_B}{dt} = k_0 - k_1 E_B$. $\qquad (4)$

Because in isotropic space $\forall i : E_{B_i} = const$, meaning $\frac{dE_B}{dt} = 0$, it follows that $k_0 = 0$.

Hence $\int \frac{dE_B}{E_B} = const - k_1 \int dt$, it follows that solution of (4) will look like

$$E_B = E_{B,0} \exp(-k_B t), \qquad (5)$$

where $E_{B,0}$ is the energy concentration in B point particle at the initial time moment; $k_1 = k_B$ is a constant presenting energy emission intensity by the point, corresponding to B body.

It follows from (5) that $m_B$ mass value concentrated in B particle is not constant:

$$m_B = \frac{E_{B,0}}{C^2} \exp(-k_B t) = m_{B,0} \exp(-k_B t) \qquad (6)$$

Consequently: $\qquad m_B = -\frac{1}{k_B} \frac{dm_B}{dt}. \qquad (7)$

Thus $m_B$ *mass of material body B can be determined by mass loss intensity of this body.* For short time intervals (5) and (6) can be substituted by:

$$E_B = E_{B,0}(1 - k_B t), \qquad (8)$$
$$m_B = m_{B,0}(1 - k_B t) \qquad (9)$$

From equations (5) - (9) we get a physical meaning of $k$ emission coefficient which is the relative change of energy/mass in a selected body in a time unit (t):

$$k = \frac{E_0 - E}{E_0 t} = \frac{m_0 - m}{m_0 t}. \qquad (10)$$

Then change of $\Delta E_B$ energy or of $\Delta m_B$ mass in B body after the expiration of $t$ time is determined by system:

$$\begin{cases} \Delta E_B = E_{B,0} - E_B(t) = E_{B,0} k_B t; \\ \Delta m_B = m_{B,0} - m_B(t) = m_{B,0} k_B t. \end{cases} \qquad (11)$$

Let us suppose that in $K_0$ coordinate system there are two motionless bodies A and B with masses of $m_A$ and $m_B$, respectively, their emission coefficients are coinciding. It is obvious that an internal observer within this system using as a reference mass or energy of one of these bodies will not be able to observe changes of these parameters for the other body because of equations (5), (6), (8) and (9) due to:

$$\begin{cases} \dfrac{E_A}{E_B} = \dfrac{E_{A,0}}{E_{B,0}}; \\ \dfrac{m_A}{m_B} = \dfrac{m_{A,0}}{m_{B,0}}. \end{cases} \quad (12)$$

## 3. Newton Laws

The aim of this section is at developing of system of models, permitting to get Newton Lows or their equivalents for a system of emitting energy bodies using basic conceptions above.

### 3.1. First Newton Law

If the space adjacent to A body is isotropic, i.e. A does not experience any external effects then the F forces resultant of $F_s$ reaction being the consequence of energy emission from each unit of s surface is determined by $F = \oint F_s ds$. «Spanning» body A into a point A it is easy to see (Figure 1) that for any reaction force vector of $\mathbf{F}(\varphi)$ corresponding to φ direction in relation to A there is a similar force in value and in opposite in direction $\mathbf{F}(-\varphi) = -\mathbf{F}(\varphi)$, i.e. resultant of reaction forces F = 0.

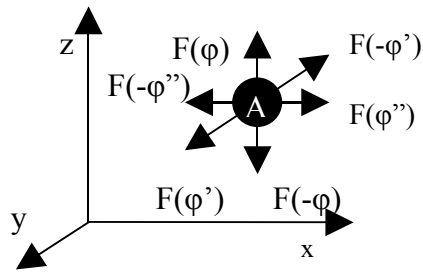

Figure 1. Reaction forces in isotropic medium.

In other words *reaction forces in the isotropic medium are unable to change the state of A material point if there are no external effects, this is found to correspond to the First Newton Law.*

### 3.2. Second Newton Law

Let A energy emitting body correspond to a certain material point A, moving with $a_A$ constant acceleration during a time interval (*t*). Keeping in mind (5), $F_A$ acceleration force is equal to:

$$F_A = a_A m_{A,0} \exp(-k_A t). \quad (13)$$

For small *t* values in accordance with (11) equation (13) is transformed to become:

$$F_A = a_A m_{A,0}(1 - k_A t). \quad (14)$$

Denoting the $a_A m_{A,0}$ product by $F_{A,0}$ variable, equation (13) which characterises a force affecting $a$ acceleration to A body can be changed to:

$$F_A = F_{A,0}\exp(-k_A t). \qquad (15)$$

It means that the force imposing a constant acceleration on A body is not constant. Nevertheless it is not possible to find the inconstancy of this force by comparing two such forces affecting the energy emitting A and B bodies which have close emission coefficients $(k_A \approx k_B = k)$, due to their ratio being constant:

$$\frac{F_A}{F_B} = \frac{F_{A,0}\exp(-kt)}{F_{B,0}\exp(-kt)} = \frac{F_{A,0}}{F_{B,0}}. \qquad (16)$$

Thus taking into account the fact that for an internal observer in the system of material points for which emission coefficients have similar values the second Newton's law is valid: *a body moving with a uniform acceleration in the K coordinate system during a certain time period (t) is under the impact of a constant force (F)*.

### 3.3. Third Newton Law

To analyse the inertia we shall substitute A material body by a sphere B with $m_B$ mass and R radius (Figure 2a). It is then easy to see that the geometry of distribution of energy emitting sphere both within the shell and outside it coincides with its distribution within a cocoon formed by A body with the same mass, with the exception of particles belonging to the sphere surface and sphere's centre.

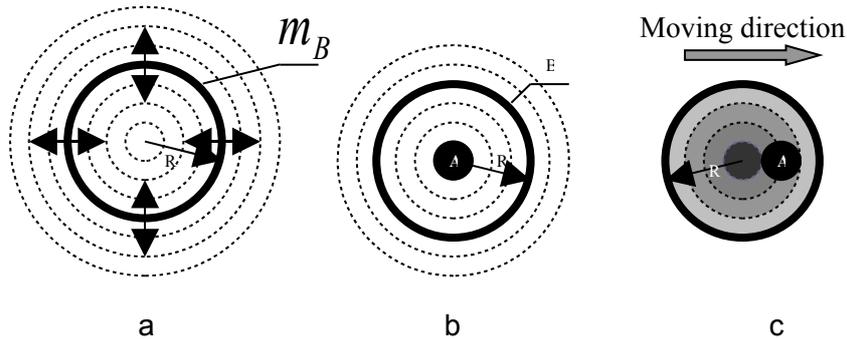

a              b              c
Figure 2. Model of reaction force initiation under accelerated motion.

Denoting density of energy distribution at distance $R_i$, $(\forall\, i \geq 1, R_i < R)$ from the centre of the sphere as $\theta(R_i)$ we can show that:

$$\begin{cases} \forall\ 0 < R_i \neq R,\ \theta(R_i) = \dfrac{k_B m_B}{4\pi R_i^2}; \\ \theta(0) = k_B m_B. \end{cases} \qquad (17)$$

It allows us to formally detach the energy cocoon surrounding motionless material point (A) from the particle proper, substituting the cocoon's layers (which are concentric hollow spheres) by the distribution of energy emitted by B sphere, the distribution being determined by system (17) provided that $m_A = m_B$, $k_A = k_B$ (Figure 2b). The accelerated motion of A material point during a short time interval ($t$) results in the cocoon's shift relative to A in the direction opposite to that of the acceleration, which in turn results in anisotropy of the space around A: and according to (17) leads to reduction of energy concentration density dispersed in front of A body in the direction of movement and a similar growth of energy dispersed in the space abandoned by this body (Figure 3c). The latter distorts the reaction force equilibrium, their resultant now differs from zero and according to the first two assumptions made in Section 2 the anisotropy leads to a higher intensity of A material point energy emission towards the acceleration direction and reduction of that acting in the opposite direction. Assuming $t$ time during which F force imposes $a$ acceleration to A particle and sphere radius (R) being $\frac{at^2}{2} < R$ the $F$ value according to [1] will be described by:

$$F = a\left(m_A + \gamma \frac{m_A m_B}{Rc^2}\right), \tag{18}$$

where "$\gamma$" is the gravitational constant and "$c$" is the velocity of light. Considering that the masses of B sphere and A material particle are equal we can express (18) as:

$$F = am_A\left(1 + \gamma \frac{m_A}{Rc^2}\right) \tag{19}$$

Then the value of $I$ impulse equals to $I = atm_A\left(1 + \gamma \frac{m_A}{Rc^2}\right).$ (20)

Let us define by $m_e$ the mass of non-balanced energy emitted with $V_e$ speed by A material body in a time unit ($t$) during A material body motion. We shall further assume $m_e$ to be a polynomial function of the mass and of its acceleration expressed as:

$$m_e = \sum_{i \geq 0} \sum_{j \geq 0} k_{i,j} m_A^i a^j, \tag{21}$$

i.e. in the first approximation if t→0, then according to our assumption 3 above we can transform (21) as $m_e = k_{0,0} + k_{1,0}m_A + k_{0,1}a + k_{1,1}am_A.$ (22)

Since coefficient values in (22) do not depend on the acceleration and mass of A body we can further analyse three combinations of these variables values:

1. $m_A = 0$, $a \neq 0$. Obviously in this case $m_e=0$, from which it follows $k_{0,0} + k_{0,1}a = 0$.
2. $m_A \neq 0$, $a = 0$. Obviously in this case $m_e=0$, from which it follows $k_{0,0} + k_{1,0}m = 0$.
3. It follows from $m_A = a = 0$ that $m_e = 0$ and we further obtain $k_{0,0} = 0$.

System:

$$\begin{cases} k_{0,0} + k_{0,1}a = 0; \\ k_{0,0} + k_{1,0}m = 0; \\ k_{0,0} = 0, \end{cases}$$

is compatible only if equalities $k_{0,0} = k_{1,0} = k_{0,1} = 0$ are valid. This corresponds to the change of (22) by:

$$m_e = k_{1,1}am_A. \tag{23}$$

Using equation (23) we shall determine the change of linear momentum (Q) caused by $m_e$ emission during $t$ time with a speed of $V_e$:   $Q = k_{1,1}am_A tV_e$   (24)

We can determine the value of $k_{1,1}$ coefficient by equating second members of equality (20) and (24): $k_{1,1} = \frac{1}{V_e}\left(1 + \gamma \frac{m_A}{Rc^2}\right).$ (25)

Now we shall determine R value with which the developed model is adequate. Since the first time derivative of A body linear momentum is equal to F force under impact of which this body is moving with a constant acceleration (*a* ) in a time interval (*t*) and reaches v speed, value of F can be presented as:

$$F = am_{0,A} \frac{1 + \frac{v^2}{c^2}}{\left[1 - \frac{v^2}{c^2}\right]^{\frac{3}{2}}}. \qquad (26)$$

Equating second members of equations (19) and (26), provided that t→0, we can determine R value:

$$R = \frac{\gamma m_A}{c^2} \left[ \frac{1 + \frac{v^2}{c^2}}{\left(1 - \frac{v^2}{c^2}\right)^{\frac{3}{2}}} - 1 \right]^{-1}. \qquad (27)$$

By substituting the second member of (27) in (25), we get: $\lim_{t \to 0} k_{1,1} = \frac{1}{V_e}$. (28)

Now by substituting the second member of (28) in (24), and of (27) – in (20), it is easy to see that within the framework of the developed model for short time intervals the values of *I* and *Q* coincide, while *the acceleration and reaction forces are equal and are acting in opposite directions, which corresponds to the third Newton's Law.* Below we shall analyse the properties of a system of two distant energy emitting bodies by means of an imaginary experiment.

### 3.4. The law of gravity

Let two emitting energy bodies (A and B) to be located at R distance from each, the interaction being manifested by energy flow from one body which reaches the layer closest to the surface of the energy cocoon formed by another body (Figure 3), this flow is added to the cocoon energy formed by this body.

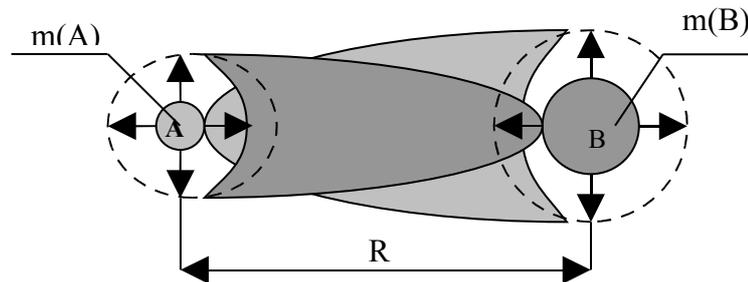

Figure 3. Remote interaction of A and B bodies, dotted lines are indicating the cocoon layers closest to the surface.

`The media around both bodies is obviously anisotropic: denoting energy emitted by D-body ($D \in \{A,B\}$) during *t* time by $E_e(D) = m_e(D) \cdot c^2$, where $m_e(D)$ - is a mass corresponding to this energy, c – velocity of light, then using the theorem of speed sums we can show that D increases kinetic energy of the cocoon of the other body in its direction by a value equal to $\Delta E_e^+(D)$:

$$\Delta E_e^+(D) = \frac{E_e(D)}{8\pi R^2 c^2} \left[ \frac{V_D + V_{\bar{D}}}{1 + \frac{V_D V_{\bar{D}}}{c^2}} \right]^2, \qquad (29)$$

and equal to $\Delta E_e^-(D)$, – in the opposite direction:

$$\Delta E_e^-(D) = \frac{E_e(D)}{8\pi R^2 c^2} \left[ \frac{V_D - V_{\bar{D}}}{1 - \frac{V_D V_{\bar{D}}}{c^2}} \right]^2, \qquad (30)$$

where $\bar{D} = \{A,B\} \setminus D$, $V_{\bar{D}}$ - the speed of $m_e(D)$.

For $t$ time the portion of the energy lost by D body that reached $\bar{D}$ nearest to the surface layer of its energy cocoon will be added to its energy $E_e^*(D) = \Delta E_e^+(D) + \Delta E_e^-(D)$ distorting the isotropy of the surrounding $\bar{D}$ space (Figure 3). As a result the balance of reactive forces is violated and $F(\bar{D})$ reactive force appears which affects $\bar{D}$ body, being directed towards D. Under the impact of $F(\bar{D})$ force $\bar{D}$ body moving towards D, during $t$ time fulfils work determined by:

$$W(\bar{D}) = F(\bar{D}) V(\bar{D}) t = \Delta E_e^+(D) - \Delta E_e^-(D), \qquad (31)$$

where $V(\bar{D})$ - $\bar{D}$ body average speed of motion in $t$ time.

If we assume that for all material points the following is valid:
- $V_D$ speed is a constant ($\forall D : V_D = const.$),
- $\forall D : k_D = const = k$,

it follows from (31), on the basis of (9), that:

$$\forall D \in \{A,B\} : F(D) = \frac{m_D k V_D^2}{2\pi R^2 V(\bar{D}) \left[1 + \frac{V_D^2}{c^2}\right]^2}. \qquad (32)$$

System (32) can be transformed to look as:

$$\forall D \in \{A,B\} : F(D) = K \frac{m_D}{V(\bar{D})}, \qquad (33)$$

where:

$$K = k \frac{V_D^2}{2\pi R^2 \left[1 + \frac{V_D^2}{c^2}\right]^2}.$$

Thus under the action of F(A) and F(B) reactive forces A and B bodies would move towards each other as if they were mutually attracted, $E_k$ kinetic energy received during $t$ time by each body considering (29), being equal respectively to:

$$\begin{cases} E_k(A) = \dfrac{m_A[V(A)]^2}{2} = \dfrac{m_{B,0}ktc^2}{4\pi R^2}; \\ E_k(B) = \dfrac{m_B[V(B)]^2}{2} = \dfrac{m_{A,0}ktc^2}{4\pi R^2}. \end{cases} \quad (34)$$

Assuming that for small $t$ we can have $m_{A,0} = m_A$, $m_{B,0} = m_B$, from ratio $\dfrac{E_k(A)}{E_k(B)}$ on the basis of (34) we obtain: $|m_A V(A)| = |m_B V(B)|$. (35)

Ignoring opposite direction of speed vectors of A and B bodies, the time derivatives of the second and first members of (35) are equal respectively to:

$$\begin{cases} \dfrac{d}{dt}[m_A V(A)] = m_A \dfrac{dV(A)}{dt} - V(A)\dfrac{k}{C^2}m_A = m_A a_A - V(A)\dfrac{k}{C^2}m_A; \\ \dfrac{d}{dt}[m_B V(B)] = m_B \dfrac{dV(B)}{dt} - V(B)\dfrac{k}{C^2}m_B = m_B a_B - V(B)\dfrac{k}{C^2}m_B, \end{cases} \quad (36)$$

where: $a_j = \dfrac{dV(j)}{dt}$ - acceleration of $j$-th solid ($j$=A, B). Since $m_j a_j = F(j)$, from (36) we obtain the equality: $F(A) - V(A)\dfrac{k}{C^2}m_A = F(B) - V(B)\dfrac{k}{C^2}m_B$. (37)

It follows from (35) that: $m_A V(A)\dfrac{k}{C^2} = m_B V(B)\dfrac{k}{C^2}$, then on the basis of (37), we get the *law of reactive forces equality*, in accordance with which the A and B bodies are moving to each other:

$$|F(A)| = |F(B)|. \quad (38)$$

Apparently, *(38) is an analogue of the law of equality of gravitational forces, impacting an isolated system of A and B bodies.*

Based on (33), equality (38) can be transformed to:

$$\dfrac{m_A}{V(B)} = \dfrac{m_B}{V(A)}. \quad (39)$$

Assuming that for short $t$ time intervals with fixed R distance any ratio of (39), being the function of A and B solids can be determined by polynomial of $f(m_A, m_B)$ type:

$$\dfrac{m_A}{V(B)} = \dfrac{m_B}{V(A)} = f(m_A, m_B) = \sum_{i \geq 0}\sum_{j \geq 0} b_{i,j} m_A^i m_B^j ,$$

or, in the first approximation and according to assumption 3 in Section 2:

$$\dfrac{m_A}{V(B)} = \dfrac{m_B}{V(A)} = f(m_A, m_B) = b_{0,0} + b_{1,0}m_A + b_{0,1}m_B + b_{1,1}m_A m_B. \quad (40)$$

Since in (40) the values of coefficients $b_{i,j}$ ($i$ =1,0; $j$=1, 0), do not depend on values of $m_A$, $m_B$ masses, below we consider three versions of the combination of these values:

1. Let $m_A = m_B = 0$. Then $f(m_A, m_B) = 0$, we get: $b_{0,0} = 0$.
2. Let $m_A = 0$, $m_B \neq 0$. Then $f(m_A, m_B) = 0$, we get: $b_{0,0} = 0$ and $b_{1,0} = 0$.
3. Let $m_A \neq 0$, $m_B = 0$. Then $f(m_A, m_B) = 0$, we get: $b_{0,0} = 0$ and $b_{0,1} = 0$.

We thus obtain:

$$\frac{m_A}{V(B)} = \frac{m_B}{V(A)} = f(m_A, m_B) = b_{1,1} m_A m_B = b m_A m_B, \tag{41}$$

where $b_{1,1} = b$. Then the following is valid: $\forall j \in \{A,B\} : V(j) = \frac{1}{bm_j}.$ (42)

Returning to (32) - (33), and on the basis of (42), we get:

$$\forall D \in \{A,B\} : F(D) = \frac{m_A m_B kb}{2\pi R^2} \left[1 + \frac{V_D^2}{c^2}\right]^{-2} = \varphi \frac{m_A m_B}{R^2}, \tag{43}$$

where: $\varphi = \frac{kb}{2\pi}\left[1 + \frac{V_D^2}{c^2}\right]^{-2}.$ (44)

It is obvious that, with the exception of the designations accuracy, expression (43) *is analogous to Newton's law on gravitation*.

### 4. Consideration of relativistic effects

It is easy to see that the value of energy ($E_e$) lost by A body in a time unit determined with consideration of relativistic effects on the basis of (16) is invariant with respect to the speed (v) that this solid has developed during $t$ time:

$$E_e = k \frac{m_0 C^2}{\sqrt{1 - \frac{v^2}{C^2}}} t \sqrt{1 - \frac{v^2}{C^2}} = E_0 k t. \tag{45}$$

We can show that with respect to unstable body moving uniformly and linearly with respect to an observer with $v$ speed and within the framework of the above developed formalism the laws of conservation of energy and impulse are valid. Thus with the consideration of relativistic effects at a certain moment in time the entire energy of such a system consists of potential energy $E_p$ and energy $E_e$, lost by this body due to its emission:

$$E = E_p + E_e = \frac{m_0 C^2}{\sqrt{1-\frac{v^2}{C^2}}}\left[1 - \frac{kt\sqrt{1-\frac{v^2}{C^2}}}{C^2}\right] + \frac{m_0}{\sqrt{1-\frac{v^2}{C^2}}} kt \sqrt{1-\frac{v^2}{C^2}} = \frac{m_0 C^2}{\sqrt{1-\frac{v^2}{C^2}}}. \tag{46}$$

Apparently $E$ is not a function of $t$ time. In a similar way we can show that an entire impulse of force ($I$) is independent of $t$ time. Here we shall assume that $I_b$ body impulse after $t$ time and $I_e$ impulse relevant to the energy lost in this time period ($t$) relevant to the energy lost due to the mass difference ($m_0 - m$) are the components of $I$:

$$I = I_b + I_e = \frac{m_0 v}{\sqrt{1-\frac{v^2}{C^2}}}\left[1 - \frac{kt\sqrt{1-\frac{v^2}{C^2}}}{C^2}\right] + \frac{m_0}{\sqrt{1-\frac{v^2}{C^2}}} \frac{vkt\sqrt{1-\frac{v^2}{C^2}}}{C^2} = \frac{m_0 v}{\sqrt{1-\frac{v^2}{C^2}}}. \tag{47}$$

### 5. Close interactions, binding energy and emission coefficients

Below we determine the values of emission coefficients (k) with respect to the nuclei of the Periodic Table elements on the basis of the liquid-drop nuclear model with the following assumptions [2]:
- all nucleons in the first approximation are spherical and all have a similar radius of $r \approx 0.8$ fm.;
- distance L, at which the interaction between nucleons is "felt", does not exceed 1 fm.;
- energy emitted by neighbouring nucleons, only one of them belonging to the surface of an atomic nucleus, in the direction of each other results in the anisotropy of energy distribution in the space around the "surface nucleons" and as a consequence of a resultant reaction force of each such a nucleon which is different from zero.

Figure 4 shows nucleons of the nuclei of deuterium and tritium as continuous circles marked either "p" for protons, or "n" for neutrons. Dotted line marks the areas within which the nucleons are interacting.

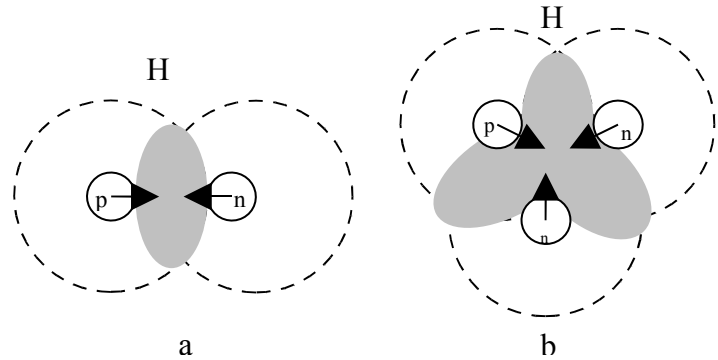

a          b

Figure 4. Resultants of reaction forces of nucleons of deuterium and tritium nuclei.

Reduction of energy emission by nucleons belonging to the nucleus in the direction of its higher concentration at distances not exceeding 1 fm. (in the Figure 4 such "contact areas" between nucleons of deuterium and tritium nuclei are shown by grey colour) on the basis of (7) results in a mass of the atom nucleus which is smaller than the total mass of "free" nucleons of that nucleus. Since the resultant of reaction forces (indicated by arrows in Figure 4) tends to facilitate approaching of each nucleon with its direct neighbours, a defect of mass is determined by the sum of "contributions" of the nucleons of this nucleus to the binding energy $\Delta E$. Assuming that in the first approximation the energy spent by j-th nucleon for the binding with its neighbours is linearly dependent on their number we shall get:

$$\begin{cases} \sum_{1 \leq j \leq A} \delta E_j = \Delta E; \\ \forall j : \delta E_j = \dfrac{\nu_j}{4\pi} \dfrac{dE_j}{dt}, \end{cases} \quad (48)$$

where: $\dfrac{dE_j}{dt}$ - energy, emitted by j-th nucleon in a time unit outside the nucleus;

$\delta E_j$ - energy, spent by j-th nucleon for binding with its neighbours;

$\nu_j$ - number of neighbours of j-th nucleon in the nucleus which form with it the areas of higher concentration of energy ($\nu_j \leq min[A-1; 4\pi]$);

A – number of nucleons in the nucleus.

Thus it is the nucleus geometry that is important in determination of binding energy, it determines the number of direct neighbours of each nucleon, the distance for which does not exceed L.

In order to estimate the value of k emission coefficient we shall investigate several geometric interpretations of the liquid-drop nuclear model, each interpretation being valid for "its own" range of A value changes. One of these models, the simplest one, allows us to estimate the emission coefficient value (k) with respect to the nuclei of light elements of the Periodic Table for which the following is valid:

$$\begin{cases} \forall j : v_j = const; \\ \forall j : 1 < A_j < 5. \end{cases}$$

Nuclei of deuterium, tritium and helium-4 are meeting these conditions. That is why for any of these nuclei the equality $\forall j : v_j = A - 1$ is valid and from (48) it follows:

$$\forall i \in \{H_1^2, H_1^3, He_2^4\}, \Delta E_i = \frac{(A_i - 1)}{4\pi} \sum_{j=1}^{A_i} \delta E_j = \frac{(A_i - 1)}{4\pi} [n_p \frac{dE_p}{dt} + (A_i - n_p) \frac{dE_n}{dt}], \quad (49)$$

where: $A_i$ - number of nucleons in i-th nucleus,

$n_p$ - number of protons in a nucleus,

$\Delta E_i$ - binding energy of the i-th nucleus in atomic units of a mass,

$\frac{dE_p}{dt}$ - energy, emitted by a proton in a time unit,

$\frac{dE_n}{dt}$ - energy emitted by a neutron in a time unit.

On the basis of (8), expression (49) can be written as:

$$\forall i \in \{H_1^2, H_1^3, He_2^4\}, \Delta E_i = \frac{931.04(A_i - 1)k_i}{4\pi} [n_p m_p + (A_i - n_p) m_n], \quad (50)$$

where: $m_p$ - proton mass in the atomic mass units,

$m_n$ - neutron mass in the atomic mass units,

931.04 is the factor of transformation of a unit of atomic mass into energy, expressed in MeV [2].

Further addition of nucleons to a nucleus, preserving its compactness, results in their "contact inequality": part of the nucleons has common "contact areas" with three others, while the other part has such common areas with a large number of neighbours. It is possible to show that this approach allows us to transform (50) for the light nuclei, i.e. for those for which $A_i \leq 50$ as follows:

$$\forall 2 < A_i \leq 50, \Delta E_i = \frac{931.04 g_i k_i}{4\pi} [n_p m_p + (A_i - n_p) m_n], \quad (51)$$

where $g_i$ is determined by:

$$g_i = \begin{cases} A_i - 1, \text{ если } A_i \leq 4; \\ [12 + 6(A_i - 4)] / A_i, \text{ если } 4 < A_i < 19; \\ [96 + 8(A_i - 19)] / A_i, \text{ если } 19 \leq A_i \leq 25; \\ [136 + 6(A_i - 24)] / A_i, \text{ если } A_i > 25. \end{cases} \quad (52)$$

It allows us to use (53) to determine $k_i$ for the nuclei meeting the constraints of (52):

$$\forall\, A_i > 2,\ k_i = \frac{0.01349 \cdot \Delta E_i}{g_i \cdot [n_p m_p + (A_i - n_p) m_n]}. \qquad (53)$$

Size of emission coefficient k changing, obtained on the basis of (53) for atomic nuclei of 80 elements of the Periodic Table and their isotopes containing from 2 to 144 nucleons is reflected by range $R_1 = [0.0139 - 0.03153671]$, average coefficient k value is equal to $k^* = 0.01945$. The use of only stable elements, their number in the sample being 55, only slightly changed the value: $k^{**} = 0.01972$. Appendix 1 presents Table of Emission Coefficients ($k_i$), its first and fourth columns determine $i$ index $i$; second and fifth are reflecting symbols of corresponding elements of the Periodic Table, third and sixth columns present $k_i$ values corresponding to these elements.

Let us divide the set of all the material points into n such sub-sets, that each $i$-th of them ($i = 1, 2, \ldots n$) would have a corresponding sub-set of emission coefficients ($K_i = \bigcup_j k_{i,j}$) which make $\forall\, 1 \leq i \leq n, \forall\, j \neq q, k_{i,j} \approx k_{i,q}$. In other words any pair of material bodies A and B whose emission coefficients are in the same sub-set, satisfy the conditions of (12), i.e. they are stable from the point of view of an internal observer. Now we shall select $K_w$ sub-set for which (54) is valid:

$$\frac{1}{|K_w|} \sum_{j=1}^{|K_w|} k_{w,j} = \max_{\leq 1 \leq i \leq n} \frac{1}{|K_i|} \sum_{j=1}^{|K_i|} k_{i,j}. \qquad (54)$$

It is easy to see that for any pair of material points (A and B) whose emission coefficients are $k_A \in K_w, k_B \notin K_w$ A body is stable as seen by an internal observer, while B body is not stable. $K_w$ sub-set found by means of (54) is of particular importance because as seen by an internal observer there are no bodies in nature whose mass would grow spontaneously. Since it is in conformity with the observed physical processes the suggested concept based on models of energy emission by all the material points at first glance is not at variance with the physical lows. Experimental validation of these models is described below.

## 6. Experiment development and its results

Experimental validation of the efficiency of the above approaches was made in two ways. First we verify their use in order to improve the quality of predicting the binding energy of stable nuclei using semi-empirical binding energy equation [2]: due to small dispersion of k coefficient, calculation precision of the binding energy for light nuclei by means of original equation can be improved by substituting the summand, which determines the volumetric nucleus energy, by the second term of (51):

$$\forall\, 2 < A_i < 50,\ \Delta E_i = \frac{931.04 \cdot g_i \cdot k^{**}}{4\pi} [n_p m_p + (A_i - n_p) m_n] + s + h + y, \qquad (55)$$

where $g_i$ is determined by (52), whereas the other components – by the original equation [2]:

$$s = -0.083 \cdot \frac{\left(\frac{A}{2} - n_p\right)^2}{A} - 627 \cdot 10^{-6} \cdot \frac{n_p(n_p - 1)}{\sqrt[3]{A}}, \tag{56}$$

$$h = 0.14\sqrt[3]{A^2} - 81 \cdot 10^{-5} \cdot n_p, \tag{57}$$

$$y = \begin{cases} 0, & \text{for odd-mass nucleus}; \\ \frac{33.517}{\sqrt[4]{A^3}} - 2, & \text{for even-even nucleus}; \\ -\frac{33.517}{\sqrt[4]{A^3}}, & \text{for odd-odd nucleus}. \end{cases} \tag{58}$$

First two columns of Table 2 (Appendix 2) present the line number and corresponding to the line shorthand notation of a stable nuclide, respectively; third, fourth and fifth columns present binding energy values obtained experimentally (3-d column), by (54) (fourth column), and on the basis of an original semi-empirical formula (5-th column). Grey colour represents lines in which the absolute value of the difference between the values of the third and fourth columns does not exceed that of the difference between the values of third and fifth columns. Thus in 70 % cases the modified formula gave better results then the original one. More accurate estimate of comparative efficiency of these two formulas would be obtained by means of a method of comparison standards [3] in which the components of the third column play the role of standard values. However its use favours (54).

The second approach aimed at verification of the range of emission coefficient values of stable or long-lived nuclei, the range being determined by means of radioactive elements. While measuring the mass of a radioactive element ($m_i$) by means of a standard mass ($m_j$) of a stable nuclei, it should be remembered that the latter also changes in time, i.e. the ratio of these masses on the basis of (9) is:

$$\frac{m_i}{m_j} = \frac{m_{0,i}}{m_{0,j}} exp[(k_j - k_i)t] \tag{59}$$

Choosing as $t$ time interval the half-value period $\tau_i$ of i-th radioactive element on the basis of (59) we get:
$$ln\,0.5 = (k_j - k_i)\tau_i. \tag{60}$$

Now we then can determine $k_j$ by $k_i$ and $\tau_i$:  $k_j = k_i + \frac{ln\,0.5}{\tau_i}. \tag{61}$

In order to reduce the dispersion of predicted $k_j$ values we have analysed during our experiment only those radioactive nuclei:
- the mass of which measured in atomic units does not change in time,
- nuclei which have only one type of decay.

Thus, elements with α – decay were excluded for the analysis, this group nuclei those with ±β –decay. The first column of Table 3 (Appendix 3) shows the numbers of lines, the second one –selected radioactive elements, the third column has their emission coefficients ($k_i$) and the fourth – corresponding $\tau_i$ half lives in seconds, while the fifth column presents their corresponding $k_j$ emission coefficients of stable nuclei estimated by means of (61). It is easy to see that above estimated $k^*$ and $k^{**}$ values appear to be within the range of ψ = [0.0136005 – 0.0198447], determined by the last column of Table 3 (Figure 5). Emission coefficients range χ = [1.67848·10$^{-2}$ - 3.153671·10$^{-2}$] corresponds to stable nuclei and to equation (54) usage.

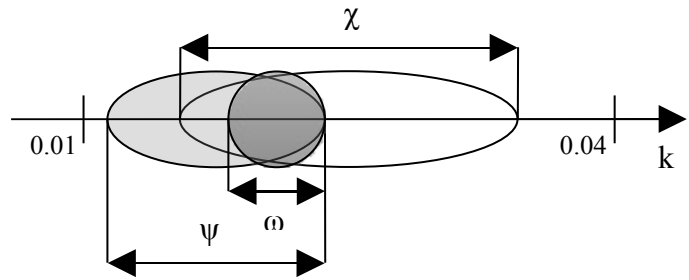

Figure 5. Emission coefficients ranges χ , ψ and ω location on the number axis.

Further contraction of this range due to the use of nuclei with only -β-decay in (61) (their corresponding lines in Table 3 are shown in grey) would give us no new result: both average values $k^*$ and $k^{**}$ lie in a new range ω = [0.01770958 - 0.0198447] (Figure 5).

## 7. Conclusion

Experimental verification of proposed above concept in the first approximation confirms its validity, while we cannot as yet consider it statistically significant because analysed were only 55 out of 285 stable and long-lived nuclei and only 14 out of more than 2000 radioactive ones. Thus it is obvious that it is reasonable to continue the research using both types of verification, based on extension of the range of nuclei used.

The other verification possibility can be connected with calculation of binding energy for nuclei using semi-empirical binding energy equation and refining model of Coulomb repulsion taking into account combinatorial analyse of different distributions of protons in a nucleus.

## 8. References

1. A. Einstein. Gibt es eine Gravitationswirkung, die der elektrodynamischen Induktionswirkung analog ist? Viertejahrschr. gerichtl. Med., 1912, Ser. 3, 44, s. 37 – 40.
2. E. Shpolsky. Atomic physics, Volume II, State Publishing House of Technical and Theoretical Literature, Moscow, 1951, pp. 742 – 744. (russ).
3. V.O. Groppen. New Solution Principles of Multi-Criteria Problems Based on Comparison Standards. www.arxiv.org/ftp/math/papers/0501/0501357.pdf, 2004.

| $i$ | $i$-th nuclide | $k_i$ | $i$ | $i$-th nuclide | $k_i$ |
|---|---|---|---|---|---|
| 1 | 2 | 3 | 4 | 5 | 6 |
| 1 | $H_1^2$ | 0.01487415 | 41 | $Cl_{17}^{35}$ | 0.01974442 |
| 2 | $H_1^3$ | 0.01891915 | 42 | $Cl_{17}^{37}$ | 0.01981747 |
| 3 | $He_2^3$ | 0.0172217 | 43 | $Ar_{18}^{36}$ | 0.01972169 |
| 4 | $He_2^4$ | 0.03153671 | 44 | $Ar_{18}^{38}$ | 0.01984931 |
| 5 | $He_2^5$ | 0.02028438 | 45 | $K_{19}^{38}$ | 0.01949287 |
| 6 | $He_2^6$ | 0.01610401 | 46 | $Ar_{18}^{40}$ | 0.0196945 |
| 7 | $Li_3^6$ | 0.01782666 | 47 | $K_{19}^{39}$ | 0.01973092 |
| 8 | $Li_3^7$ | 0.01749358 | 48 | $K_{19}^{40}$ | 0.01962977 |
| 9 | $Be_4^7$ | 0.01664449 | 49 | $Ca_{20}^{40}$ | 0.01971896 |
| 10 | $Be_4^8$ | 0.02086939 | 50 | $Ar_{18}^{41}$ | 0.01397156 |
| 11 | $Be_4^9$ | 0.0185187 | 51 | $Ca_{20}^{42}$ | 0.01978483 |
| 12 | $B_5^{10}$ | 0.01804123 | 52 | $Ca_{20}^{43}$ | 0.01969448 |
| 13 | $B_5^{11}$ | 0.01887152 | 53 | $Ca_{20}^{45}$ | 0.01984475 |
| 14 | $C_6^{12}$ | 0.02054377 | 54 | $Sc_{21}^{45}$ | 0.01985941 |
| 15 | $C_6^{13}$ | 0.01967787 | 55 | $Ti_{22}^{46}$ | 0.01982922 |
| 16 | $N_7^{14}$ | 0.01944138 | 56 | $Ti_{22}^{47}$ | 0.01986394 |
| 17 | $N_7^{15}$ | 0.01969392 | 57 | $Ti_{22}^{48}$ | 0.01990625 |
| 18 | $O_8^{16}$ | 0.0203198 | 58 | $Ti_{22}^{49}$ | 0.0198121 |
| 19 | $O_8^{17}$ | 0.01958026 | 59 | $Ti_{22}^{50}$ | 0.01989175 |
| 20 | $O_8^{18}$ | 0.01947679 | 60 | $V_{23}^{51}$ | 0.02001475 |
| 21 | $F_9^{18}$ | 0.0190572 | 61 | $Cr_{24}^{52}$ | 0.02002203 |
| 22 | $F_9^{19}$ | 0.01937953 | 62 | $Cr_{24}^{53}$ | 0.0199929 |
| 23 | $Ne_{10}^{20}$ | 0.02065861 | 63 | $Fe_{26}^{54}$ | 0.01986263 |
| 24 | $Ne_{10}^{21}$ | 0.0199899 | 64 | $Fe_{26}^{55}$ | 0.01986871 |
| 25 | $Ne_{10}^{22}$ | 0.01981272 | 65 | $Fe_{26}^{56}$ | 0.01982115 |
| 26 | $Si_{14}^{31}$ | 0.01953335 | 66 | $Fe_{26}^{57}$ | 0.01979011 |
| 27 | $Na_{11}^{22}$ | 0,01946794 | 67 | $Co_{27}^{59}$ | 0.02028432 |
| 28 | $Na_{11}^{23}$ | 0.0194818 | 68 | $Xe_{54}^{130}$ | 0.01900044 |
| 29 | $Mg_{12}^{24}$ | 0.01942585 | 69 | $Sm_{62}^{144}$ | 0.01868499 |
| 30 | $Mg_{12}^{25}$ | 0.01982245 | 70 | $U_{92}^{238}$ | 0.0167848 |
| 31 | $Al_{13}^{29}$ | 0.01943383 | 71 | $C_6^{11}$ | 0.01807981 |
| 32 | $Al_{13}^{27}$ | 0.01953686 | 72 | $N_7^{13}$ | 0,01896535 |

| 1 | 2 | 3 | 4 | 5 | 6 |
|---|---|---|---|---|---|
| 33 | $Si_{14}^{28}$ | 0.01977223 | 73 | $O_8^{15}$ | 0,01910168 |
| 34 | $Si_{14}^{29}$ | 0.019577 | 74 | $P_{15}^{32}$ | 0.01968318 |
| 35 | $P_{15}^{30}$ | 0.01943093 | 75 | $Cl_{17}^{36}$ | 0.01968651 |
| 36 | $P_{15}^{31}$ | 0.01975465 | 76 | $Cl_{17}^{39}$ | 0.01956023 |
| 37 | $S_{16}^{32}$ | 0.01975493 | 77 | $Ti_{22}^{45}$ | 0.0207741 |
| 38 | $S_{16}^{33}$ | 0.01977075 | 78 | $Ga_{31}^{69}$ | 0.0198312 |
| 39 | $S_{16}^{35}$ | 0.01977043 | 79 | $Li_3^8$ | 0.01520925 |
| 40 | $Cl_{17}^{34}$ | 0.01967553 | 80 | $Al_{13}^{26}$ | 0.01924255 |

Appendix 2

| i | i-th nuclide | $\Delta E_{exp}$ | $\Delta E_i$ | $\Delta E_V$ |
|---|---|---|---|---|
| 1 | 2 | 3 | 4 | 5 |
| 1 | $H_1^2$ | 2.2241 | -16,77427 | -16,31982 |
| 2 | $He_2^3$ | 7,7243 | 9,082245 | 2,157451 |
| 3 | $He_2^4$ | 28,2937 | 27,8039 | 30,17239 |
| 4 | $He_2^5$ | 27,3 | 26,81005 | 22,14478 |
| 5 | $Li_3^6$ | 31,987 | 26,92148 | 25,56065 |
| 6 | $Li_3^7$ | 39,239 | 44,51587 | 40,6302 |
| 7 | $Be_4^9$ | 58,153 | 62,21107 | 59,20402 |
| 8 | $B_5^{10}$ | 64,744 | 65,09291 | 63,43464 |
| 9 | $B_5^{11}$ | 76,192 | 79,89884 | 77,90697 |
| 10 | $C_6^{12}$ | 92,156 | 91,93594 | 92,88525 |
| 11 | $C_6^{13}$ | 97,102 | 97,58096 | 96,70341 |
| 12 | $N_7^{14}$ | 104,653 | 101,7856 | 101,2442 |
| 13 | $N_7^{15}$ | 115,485 | 115,2586 | 115,5511 |
| 14 | $O_8^{16}$ | 127,612 | 126,2815 | 128,6358 |
| 15 | $O_8^{17}$ | 131,754 | 132,9324 | 134,4133 |
| 16 | $O_8^{18}$ | 139,798 | 143,6057 | 145,9268 |
| 17 | $F_9^{19}$ | 147,79 | 150,6031 | 153,2602 |
| 18 | $Ne_{10}^{20}$ | 160,63 | 155,1088 | 164,9919 |
| 19 | $Ne_{10}^{21}$ | 167,39 | 165,3371 | 172,0682 |
| 20 | $Ne_{10}^{22}$ | 177,76 | 178,4066 | 184,0856 |
| 21 | $Na_{11}^{23}$ | 186,44 | 188,8704 | 190,8181 |
| 22 | $Mg_{12}^{24}$ | 197,52 | 201,7219 | 201,4773 |
| 23 | $Mg_{12}^{25}$ | 204,52 | 209,4676 | 209,4948 |
| 24 | $Al_{13}^{27}$ | 224,944 | 227,1289 | 228,0859 |

| 1 | 2 | 3 | 4 | 5 |
|---|---|---|---|---|
| 25 | $Si_{14}^{28}$ | 236,52 | 236,7066 | 237,8406 |
| 26 | $Si_{14}^{29}$ | 242,97 | 244,7884 | 246,5813 |
| 26 | $Si_{14}^{29}$ | 242,97 | 244,7884 | 246,5813 |
| 27 | $P_{15}^{31}$ | 262,898 | 262,898 | 264,9727 |
| 28 | $S_{16}^{32}$ | 271,76 | 271,7596 | 273,9333 |
| 29 | $S_{16}^{33}$ | 280,85 | 280,0886 | 283,2533 |
| 30 | $Cl_{17}^{35}$ | 298,19 | 297,7577 | 301,4175 |
| 31 | $Cl_{17}^{37}$ | 317,08 | 315,4211 | 321,7042 |
| 32 | $Ar_{18}^{36}$ | 306,69 | 306,8586 | 309,6619 |
| 33 | $Ar_{18}^{38}$ | 326,49 | 324,4332 | 332,5948 |
| 34 | $Ar_{18}^{40}$ | 341.62 | 343.4461 | 351.1033 |
| 35 | $K_{19}^{38}$ | 320,62 | 322,0409 | 322,9264 |
| 36 | $K_{19}^{39}$ | 333,39 | 333,0637 | 337,3782 |
| 37 | $Ca_{20}^{40}$ | 342,03 | 341,9898 | 344,9641 |
| 38 | $Ca_{20}^{42}$ | 360,93 | 359,5778 | 368,4792 |
| 39 | $Ca_{20}^{43}$ | 368,12 | 368,3764 | 376,6555 |
| 40 | $Sc_{21}^{45}$ | 389,02 | 386,026 | 394,6564 |
| 41 | $Ti_{22}^{45}$ | 406.93 | 387.635 | 390.3501 |
| 42 | $Ti_{22}^{46}$ | 397,32 | 394,7425 | 403,818 |
| 43 | $Ti_{22}^{47}$ | 406,93 | 403,6745 | 412,5005 |
| 44 | $Ti_{22}^{48}$ | 416,73 | 412,3411 | 424,0958 |
| 45 | $Ti_{22}^{49}$ | 423,65 | 421,3308 | 431,1398 |
| 46 | $Ti_{22}^{50}$ | 434,28 | 429,99403 | 441,0839 |
| 47 | $V_{23}^{51}$ | 445,94 | 438,9781 | 449,2482 |
| 48 | $Cr_{24}^{52}$ | 455,08 | 447,528 | 459,5895 |
| 49 | $Cr_{24}^{53}$ | 463,39 | 456,6244 | 467,1747 |
| 50 | $Fe_{26}^{54}$ | 469,27 | 465,1136 | 472,8117 |
| 51 | $Fe_{26}^{56}$ | 486,08 | 482,7249 | 494,4324 |
| 52 | $Fe_{26}^{57}$ | 494,2 | 491,9146 | 502,4948 |
| 53 | $Co_{27}^{59}$ | 524,74 | 509,5585 | 519,8936 |
| 54 | $Xe_{54}^{130}$ | 1195,87 | 1257,869 | 1189,607 |
| 55 | $Sm_{62}^{144}$ | 1780 | 2086,011 | 1790,972 |

| i | $i$-th nuclide | $k_i \left(\dfrac{1}{sec}\right)$ | $\tau_i$ (sec) | $k_j \left(\dfrac{1}{sec}\right)$ |
|---|---|---|---|---|
| 1 | 2 | 3 | 4 | 5 |
| 1 | $H_1^3$ | 0.01891915 | $3.815856 \cdot 10^8$ | 0.01891915 |
| 2 | $Ca_{20}^{45}$ | 0.01984475 | $1.31328 \cdot 10^7$ | 0.0198447 |
| 3 | $C_6^{11}$ | 0.01807981 | 1340.04 | 0.01756255 |
| 4 | $N_7^{13}$ | 0.01896535 | 597.9 | 0.01780605 |
| 5 | $O_8^{15}$ | 0.01910168 | 126 | 0.01360051 |
| 6 | $F_9^{18}$ | 0.0190572 | 6420 | 0.01894924 |
| 7 | $Al_{13}^{26}$ | 0.01924255 | $2.26 \cdot 10^{13}$ | 0.01924255 |
| 8 | $Al_{13}^{29}$ | 0.01943383 | 402 | 0.01770958 |
| 9 | $P_{15}^{30}$ | 0.01943093 | 153 | 0.01490055 |
| 10 | $Si_{14}^{31}$ | 0.0195335 | 9438 | 0.01945991 |
| 11 | $P_{15}^{32}$ | 0.01968318 | 1232237 | 0.01968262 |
| 12 | $S_{16}^{35}$ | 0.01977043 | 7560864 | 0.01977034 |
| 13 | $F_9^{18}$ | 0.0190572 | 6420 | 0.01894924 |
| 14 | $Cl_{17}^{39}$ | 0.01956023 | 3600 | 0.01936769 |